\documentclass[12pt]{article}
\usepackage{amsmath,amsthm,amssymb}
\usepackage{mathrsfs}
\usepackage[titletoc]{appendix}
\usepackage{bbm}
\usepackage{float}
\usepackage{tikz}
\begin{document}
\input amssym.def
\newcommand{\singlespace}{
    \renewcommand{\baselinestretch}{1}
\large\normalsize}
\newcommand{\doublespace}{
   \renewcommand{\baselinestretch}{1.2}
   \large\normalsize}
\renewcommand{\theequation}{\thesection.\arabic{equation}}

\setcounter{equation}{0}
\def \ten#1{_{{}_{\scriptstyle#1}}}
\def \Z{\Bbb Z}
\def \C{\Bbb C}
\def \R{\Bbb R}
\def \Q{\Bbb Q}
\def \N{\Bbb N}
\def \F{\Bbb F}
\def \L{\Bbb L}
\def \l{\lambda}
\def \V{V^{\natural}}
\def \wt{{\rm wt}}
\def \tr{{\rm tr}}
\def \Res{{\rm Res}}
\def \End{{\rm End}}
\def \Aut{{\rm Aut}}
\def \mod{{\rm mod}}
\def \Hom{{\rm Hom}}
\def \im{{\rm im}}
\def \<{\langle}
\def \>{\rangle}
\def \w{\omega}
\def \c{{\tilde{c}}}
\def \o{\omega}
\def \t{\tau }
\def \ch{{\rm ch}}
\def \a{\alpha }
\def \b{\beta}
\def \e{\epsilon }
\def \la{\lambda }
\def \om{\omega }
\def \O{\Omega}
\def \qed{\mbox{ $\square$}}
\def \pf{\noindent {\bf Proof: \,}}
\def \voa{vertex operator algebra\ }
\def \voas{vertex operator algebras\ }
\def \p{\partial}
\def \1{{\bf 1}}
\def \ll{{\tilde{\lambda}}}
\def \H{{\bf H}}
\def \F{{\Bbb F}}
\def \h{{\frak h}}
\def \g{{\frak g}}
\def \rank{{\rm rank}}
\def \({{\rm (}}
\def \){{\rm )}}
\def \Y {\mathcal{Y}}
\def \I {\mathcal{I}}
\def \A {\mathcal{A}}
\def \B {\mathcal {B}}
\def \Cc {\mathcal {C}}
\def \H {\mathcal{H}}
\def \M {\mathcal{M}}
\def \V {\mathcal{V}}
\def \O{{\bf O}}
\def \AA{{\bf A}}
\def \1{{\bf 1}}
\def\Ve{V^{0}}
\def\ha{\frac{1}{2}}
\def\se{\frac{1}{16}}
\singlespace
%\doublespace
\newtheorem{thm}{Theorem}[section]
\newtheorem{prop}[thm]{Proposition}
\newtheorem{lem}[thm]{Lemma}
\newtheorem{cor}[thm]{Corollary}
\newtheorem{rem}[thm]{Remark}
\newtheorem{con}[thm]{Conjecture}
\newtheorem*{CPM}{Theorem}
\newtheorem{definition}[thm]{Definition}
\newtheorem{de}[thm]{Definition}

\begin{center}
{\Large {\bf  Vertex operator algebras associated to the Virasoro algebra over an arbitrary field}} \\

\vspace{0.5cm} Chongying Dong\footnote{Supported by a NSF grant.}
\\
Department of Mathematics\\ University of
California\\ Santa Cruz, CA 95064 \\
Li Ren\footnote{Supported in part by China Postdoctor grant 2012M521688.}\\
 School of Mathematics \\ Sichuan University\\
Chengdu 610064, China\\
\end{center}
\hspace{1.5 cm}

\begin{abstract}
The vertex operator algebras and modules associated to the highest weight modules for the Virasoro algebra over an arbitrary field $\F$ with $\ch\F\ne 2$ are studied. The irreducible modules of vertex operator algebra $L(\ha,0)_\F$ are classified. The rationality of $L(\ha,0)_\F$ is established if $\ch\F\ne 7.$
\end{abstract}

\section{Introduction}

This paper is devoted to the study of vertex operator algebras over an arbitrary field $\F$. In particular, the vertex operator algebras associated to the highest weight modules for the Virasoro algebra $Vir_\F$ are investigated.  In the case $\F=\C$ is the field of complex numbers, the Virasoro vertex operator algebras have been studied extensively in \cite{FZ}, \cite{DMZ}, \cite{W}.

The most interesting irreducible highest weight modules $L(c,h)$ for the Virasoro algebra are the unitary ones
when $\F=\C.$ It is well known from \cite{FQS} and \cite{GKO} that if $c<1,$  $L(c,h)$ is unitary if and only if $c_m=1-\frac{6}{m(m+1)}$ for $m=2,3,\cdots$
and
$$h=h_{r,s}^{m}=
 \frac {[(m+1)r-ms]^2 -1 } {4m(m+1)} \ \ \
 (r,s \in {\N} , 1\le s \le r \le m-1 ).$$
 Moreover, $L(c_m,0)$ is a rational vertex operator algebra and $L(c_m,h_{r,s}^m)$ are exactly the irreducible $L(c_m,0)$-modules (cf. \cite{DMZ}, \cite{W}).  If $m=3,$ we have rational vertex operator algebra $L(\ha,0)$ and its irreducible modules $L(\ha,h)$ with $h=0,\ha,\se.$

  The importance of vertex operator algebra $L(\ha,0)$ comes from a discovery \cite{DMZ} that there are $48$ commutative Virasoro vectors with central charge $\ha$ in the moonshine vertex operator algebra $V^{\natural}$ \cite{FLM}. That is, $L(\ha,0)^{\otimes 48}$ is a vertex operator subalgebra of $V^{\natural}.$ This discovery leads to the theory of code and framed vertex operator algebras as presented in \cite{M1}-\cite{M2} and \cite{DGH}, a proof \cite{D} of conjecture of Frenkel-Lepowsky-Meurman \cite{FLM} that $V^{\natural}$ is holomorphic,
  a different construction of the $V^{\natural}$ \cite{M3} and proofs \cite{DGL}, \cite{LY} of weak versions of
  uniqueness of the moonshine vertex operator algebras conjectured in \cite{FLM}. See \cite{LSY}, \cite{LS}, \cite{LYY1}, \cite{LYY2}  for the other developments along this line.

  In this paper we study $L(\ha,0)_\F$ over an arbitrary field $\F$ with $\ch \F\ne 2$ where $L(c,h)_\F$ is the irreducible highest weight module for $Vir_\F$ with the central charge $c$ and highest weight $h.$
  It is proved that $L(\ha,0)_\F$
  is a rational vertex operator algebra with three inequivalent modules $L(\ha,h)_\F$ with $h=0,\ha,\se$ if $\ch\F\ne 7.$ If $\ch\F=7,$ $L(\ha,0)_\F$ has only two inequivalent irreducible modules $L(\ha,0)_\F$ and
  $L(\ha,\ha)_\F.$ Although we believe that $L(\ha,0)_\F$ is also rational in this case but we cannot find a proof.

 Our proof uses the fermionic realization of $L(\ha,h)_\F$ with $h=0,\ha,\se$ (cf. \cite{KR}) and  the singular vectors of the Verma module $V(\ha,h)_\F$ (cf. \cite{FF}). Let $H_\F=\F a$ be one dimensional space
 with a bilinear form such that $(a,a)=1.$ One can construct a vertex operator superalgebra $V(H_\F,\Z+\ha)$
 and its $\tau$-twisted module $V(H_\F,\Z)$ where $\tau$ is the canonical automorphism of $V(H_\F,\Z+\ha)$ of order
 $2$ coming from the structure of vertex operator superalgebra. If $\ch \F\ne 7$ then
 the even part $V(H_\F,\Z+\ha)_{\bar 0}$ of $V(H_\F,\Z+\ha)$ is isomorphic to $L(\ha,0)_\F,$ the old part $V(H_\F,\Z+\ha)_{\bar 1}$ is isomorphic to $L(\ha,\ha)_\F,$ and both even and old parts of $V(H_\F,\Z)$ are isomorphic to $L(\ha,\se)_\F.$ These isomorphisms are easy in the case  $\F=\C$ as $V(H_\F,\Z+\ha)$ and  $V(H_\F,\Z)$ are unitary $Vir_\C$-modules. But for arbitrary $\F,$ a direct proof of irreducibility of $V(H_\F,\Z+\ha)_{\bar s }$ and $V(H_\F,\Z)_{\bar s}$ is necessary. The realization of $L(\ha,h)_\F$ enables us to determine the singular vectors of $V(\ha,h)_\F$ explicitly
  as $L(\ha,h)_\F$ and $L(\ha,h)$ have the same graded dimensions. The explicit expression of these singular vectors
from \cite{FF} are then used to classify the irreducible $L(\ha,0)_\F$-modules which are exactly $L(\ha,h)_\F$ with
$h=0,\ha,\se$ and to prove the rationality of $L(\ha,0)_\F.$ In the case of $\F=\C$, the classification of irreducible modules for $L(\ha,0)$ was achieved with the help of Zhu's algebra (cf. \cite{DMZ}, \cite{W}).

The case $\ch\F=7$ is exceptional  as $\ha=\se.$ While $V(\ha,0)_\F$ has exactly two singular vectors of degrees
$1$ and $6$ if $\ch\F\ne 7,$ $V(\ha,0)_\F$ has a singular vector of degree $4$ when $\ch\F=7.$ This means
$L(\ha,0)_\F$ for $\ch\F= 7$ is totally different from $L(\ha,0)_\F$ for $\ch\F\ne 7.$ Although   $V(H_\F,\Z+\ha)_{\bar s}$ and  $V(H_\F,\Z)_{\bar s}$ are still completely reducible modules for $Vir_\F$ with $\ch\F=7,$ they are not irreducible anymore. Using the singular vector of degree $4$ we can also classify the irreducible $L(\ha,0)_\F$-modules which are
$L(\ha,0)_\F$ and $L(\ha,\ha)_\F.$ Unfortunately, we cannot find all singular vectors of $V(\ha,h)_\F$
for $h=0,\ha.$

There is a very interesting phenomenon on $L(\ha,h)_\F.$ During the proof of our main results, we find an $R$-form
$L(\ha,h)_R$ of $L(\ha,h)$ such that $L(\ha,h)_\F=\F\otimes_R L(\ha,h)_R$ if $\ch\F\ne 7$ where $R=\Z[\ha]$ (see Section 2 for the definition of $R$-forms for vertex operator algebras and modules over $\C$). But this is not true if $\ch\F=7.$
This makes us suspect  that if $V$ is a rational vertex operator algebra with an $R$-form I, then $\F\otimes_R I$ is a simple,
rational vertex operator algebra except for finitely many prime $p$ with $\ch\F=p.$  But we do not have any idea on how to approach to this problem in general.

We should mention that this paper, in fact, is the first paper dealing with modular vertex operator algebras and their modules. It is expected that study of $L(\ha,0)_\F$ will help people to investigate modular vertex operator algebras in general.

The paper is organized as follows: In Section 2, we define vertex operator superalgebras and their twisted modules
over an arbitrary field $\F.$ There is an extra assumption in the definition of vertex operator superalgebra:
if $u,v$ have degrees (which are weight in the case $\ch\F=0$) $s,t,$ respectively, then $u_nv$  has degree $s+t-n-1.$ In the case that $\ch\F=0,$ this follows from the other axioms of the definition. A similar assumption is also required in the definition of ordinary modules. In Section 3, we study vertex operator algebras associated to the highest weight modules for the Virasoro algebra $Vir_\F.$  Section 4 is devoted to the study
vertex operator algebra $L(\ha,0)_\F$ with $\ch \F\ne 7.$ In particular, the irreducible modules of $L(\ha,0)_\F$ are classified and rationality is obtained. We classify the irreducible modules of $L(\ha,0)_\F$ with $\ch\F=7$ in Section $5.$

\section{Basics}
\setcounter{equation}{0}
In this section we first discuss the basics of  vertex operator superalgebras and their twisted modules (cf. \cite{B}, \cite{FLM}, \cite{FFR},  \cite{DL}, \cite{X}, \cite{LL}) over an integral domain $\mathbb{D}$ with $\ch \mathbb{D}\ne 2.$
The typical examples of $\mathbb{D}$ in this paper are fields and $\Z[\ha].$ This general definition gives us flexibility in dealing with integral forms and modular vertex operator algebras. We also define the vertex operator superalgebra $R$-forms for an  integral subdomain $R$ of $\mathbb{D}.$

We first define a super $\mathbb{D}$-module as a $\Bbb Z_{2}$-graded free $\mathbb{D}$-module
 $V=V_{\bar{0}}\oplus V_{\bar{1}}$ such that both $V_{\bar{0}}$ and $V_{\bar{1}}$ are free $\mathbb{D}$-submodules.
 As usual  we let $\tilde{v}$ be $0$ if $v\in V_{\bar{0}}$, and $1$ if  $v\in V_{\bar{1}}$.

\begin{de} {\rm A  vertex operator superalgebra $V=(V,Y,1,\omega)$ over $\mathbb{D}$ is a
$\frac{1}{2}\Bbb Z$-graded super $\mathbb{D}$-module
$$V=\bigoplus_{n\in{ \frac{1}{2}\Bbb Z}}V_n= V_{\bar{0}}\oplus V_{\bar{1}}$$
with  $V_{\bar{0}}=\sum_{n\in\Z}V_n$ and
$V_{\bar{1}}=\sum_{n\in\frac{1}{2}+\Z}V_n$
satisfying $\dim V_{n}< \infty$ for all $n$ and $V_m=0$ if $m$ is sufficiently
small.  $V$ is   equipped with a linear map
\begin{align*}
& V \to (\mbox{End}\,V)[[z,z^{-1}]] ,\\
& v\mapsto Y(v,z)=\sum_{n\in{\Z}}v_nz^{-n-1}\ \ \ \  (v_n\in
\mbox{End}\,V)\nonumber
\end{align*}
and with two distinguished vectors ${\bf 1}\in V_0,$ $\omega\in
V_2$ satisfying the following conditions for $u, v \in V,$ and $m,n\in\Z,$ $s,t\in \ha\Z:$
\begin{align*}
&  u_nv\in V_{s+t-n-1} \ \ \ {\rm for} \ u\in V_s, v\in V_t \  {\rm and }\ u_nv=0\ \ \ \ \ {\rm for}\ \  n\ \ {\rm sufficiently\ large};  \\
& Y({\bf 1},z)=Id_{V};  \\
& Y(v,z){\bf 1}\in V[[z]]\ \ \ {\rm and}\ \ \ \lim_{z\to
0}Y(v,z){\bf 1}=v;\\
& [L(m),L(n)]=(m-n)L(m+n)+\frac{1}{12}(m^3-m)\delta_{m+n,0}c ;\\
& \frac{d}{dz}Y(v,z)=Y(L(-1)v,z);\\
& L(0)|_{V_n}=n
\end{align*}
where $L(m)=\o_{ m+1}, $ that is, $Y(\o,z)=\sum_{n\in\Z}L(n)z^{-n-2};$
and the {\em Jacobi identity} holds:
$$z^{-1}_0\delta\left(\frac{z_1-z_2}{z_0}\right)
Y(u,z_1)Y(v,z_2)-(-1)^{\tilde{u}\tilde {v}}z^{-1}_0\delta\left(\frac{z_2-z_1}{-z_0}\right)
Y(v,z_2)Y(u,z_1)$$
$$=z_2^{-1}\delta\left(\frac{z_1-z_0}{z_2}\right)Y(Y(u,z_0)v,z_2).$$}
\end{de}

If $v\in V_s$ we will call $s$ the degree of $v$ and write $\deg v=s.$ If we regard $s$ as a number in the field $\F$ we call $s$ the weight of $v$

In the case $V_{\bar 1}=0,$ we have a vertex operator algebra $V$ over $\mathbb{D}.$  If $\mathbb{D}=\C,$ the assumption $u_nv\in V_{s+t-n-1}$ in the definition
is a consequence of the other axioms. If $\mathbb{D}=\F$ is a field of finite characteristic $p,$ the $L(0)$ has only $p$ distinct eigenvalues on $V$ and the assumption is necessary.

An automorphism $g$ of a vertex operator superalgebra $V$ is a $\mathbb{D}$-module automorphism of $V$ such that $gY(u,z)v=Y(gu,z)gv$ for all $u,v\in V,$ $g\1=\1,$ $g\omega=\omega$ and $gV_n=V_n$ for all $n\in \ha \Z.$ It is clear that any automorphism preserves
$V_{\bar 0}$ and $V_{\bar 1}.$ There is a special automorphism  $\tau$  such that $\tau|V_{\bar 0}=1$
and $ \tau|V_{\bar 1}=-1.$ We see that $\tau$ commutes with any automorphism.

Fix $g\in \Aut(V)$ of order $T.$  We assume that $\frac{1}{T}\in \mathbb{D}$ and $\mathbb{D}$ contains a primitive $T$-th root of unity $\eta.$ Then $V$ decomposes into eigenspaces of $g$:
$$V=\oplus_{r\in \Z/T\Z}V^{r}$$
where $V^r=\{v\in V|gv=\eta^rv\}.$

\begin{de} {\rm A weak $g$-twisted $V$-module $M$ is a free $\mathbb{D}$-module  equipped
with a linear map
$$\begin{array}{l}
V\to (\End\,M)[[z^{1/T}, z^{-1/T}]]\\
v\mapsto\displaystyle{ Y_M(v,z)=\sum_{n\in\frac{1}{T}\Z}v_nz^{-n-1}\ \ \ (v_n\in
\End\,M)}
\end{array}$$
which satisfies that for all $0\leq r\leq T-1,$ $u\in V^r,$ $v\in V,$
$w\in M,$
\begin{eqnarray*}\label{g2.11}
& &Y_M(u,z)=\sum_{n\in \frac{r}{T}+\Z}u_nz^{-n-1} \label{1/2} ;\\
& &u_lw=0 \ \ \  				
\mbox{for}\ \ \ l\gg 0\label{vlw0};\\
& &Y_M(\1,z)=Id_{M};\label{vacuum}
\end{eqnarray*}
 \begin{equation*}\label{2.14}
\begin{array}{c}
\displaystyle{z^{-1}_0\delta\left(\frac{z_1-z_2}{z_0}\right)
Y_M(u,z_1)Y_M(v,z_2)-(-1)^{\tilde{u}\tilde{v}}z^{-1}_0\delta\left(\frac{z_2-z_1}{-z_0}\right)
Y_M(v,z_2)Y_M(u,z_1)}\\
\displaystyle{=z_2^{-1}\left(\frac{z_1-z_0}{z_2}\right)^{-r/T}
\delta\left(\frac{z_1-z_0}{z_2}\right)
Y_M(Y(u,z_0)v,z_2)}.
\end{array}
\end{equation*}}
\end{de}

Using the standard argument (cf. \cite{DL}) one can prove
that the twisted Jacobi identity is equivalent to the following
associativity formula
\begin{eqnarray}\label{ea2.15}
(z_{0}+z_{2})^{k+\frac{r}{T}}Y_{M}(u,z_{0}+z_{2})Y_{M}(v,z_{2})w
=(z_{2}+z_{0})^{k+\frac{r}{T}}Y_M(Y(u,z_0)v,z_2)w
\end{eqnarray}
where $w\in M$ and $k\in\Bbb Z_{+}$ such that $z^{k+\frac{r}{T}}Y_{M}(u,z)w$ involves only nonnegative integral powers of $z,$ and
commutator relation
\begin{eqnarray}\label{g2.16}
& &\ \ \ \  [Y_{M}(u,z_{1}),Y_{M}(v,z_{2})]\nonumber\\
& &=\Res_{z_{0}}z_2^{-1}\left(\frac{z_1-z_0}{z_2}\right)^{-r/T}
\delta\left(\frac{z_1-z_0}{z_2}\right)Y_M(Y(u,z_0)v,z_2)\label{ec}
\end{eqnarray}
whose component forms is given by
\begin{eqnarray}\label{g2.17}
[u_{s},v_{t}]=\sum_{i=0}^{\infty}
\left(\begin{array}{c}m\\i\end{array}\right)(u_{i}v)_{s+t-i}
\end{eqnarray}
for all $u,v\in V$ and $s,t\in\frac{1}{T}\Z.$  One can  also deduce the usual Virasoro algebra axioms:
\begin{equation*}\label{g2.18}
[L(m),L(n)]=(m-n)L(m+n)+\frac{1}{12}(m^3-m)\delta_{m+n,0}c,
\end{equation*}
\begin{equation*}\label{g2.19}
\frac{d}{dz}Y_M(v,z)=Y_M(L(-1)v,z)
\end{equation*}
(cf. \cite{DLM1}, \cite{DLM2}) for $m,n\in\Z$  where $Y_M(\o,z)=\sum_{n\in\Z}L(n)z^{-n-2}.$

The homomorphism and isomorphism of weak twisted modules are defined in an
obvious way.

\begin{de} An {\em admissible} $g$-twisted $V$-module
is a  weak $g$-twisted $V$-module $M$ which carries a
$\frac{1}{T}{\Z}_{+}$-grading
\begin{equation}\label{g2.22}
M=\oplus_{n\in\frac{1}{T}\Z_+}M(n)
\end{equation}
satisfying
\begin{eqnarray}\label{g2.23}
v_{m}M(n)\subseteq M(n+s-m-1)
\end{eqnarray}
for homogeneous $v\in V_s.$ We will call $n$  the degree of $w\in M(n).$
\end{de}

\begin{de} Let $\mathbb{D}=\F$ is a field. Vertex operator algebra $V$ is called $g$-rational if any admissible $g$-twisted $V$-module is completely reducible. $V$ is rational if $V$ is $1$-rational.
\end{de}

\begin{de} A $\C$-gradable $g$-twisted $V$-module $M$ is a weak $g$-twisted $V$-module graded by $\C$
$$M=\oplus_{\lambda\in \C} M_{\lambda}$$
such that $M_{\lambda}$ is finite dimensional for all $\lambda,$ $M_{\lambda+m}=0$ for fixed $\lambda$ and sufficiently large $m\in\frac{1}{T}\Z$ and $v_{m}M_{\lambda}\subseteq M_{\lambda+s-m-1}$ for $u\in V_s$
and $m\in \frac{1}{T}\Z.$
\end{de}

A $\C$-gradable $g$-twisted module defined here for any field $\F$ essentially is the ordinary $g$-twisted
module when $\F=\C.$ But we cannot assume here that each $M_{\lambda}$ is an eigenspace for $L(0)$ with eigenvalue $\lambda$ as $\lambda$ does not lie in $\F.$ It is also clear from the definition of admissible and
$\C$-gradable $g$-twisted modules that a uniform grading shift gives an isomorphic module of the same type.
One can also verify that a $\C$-gradable $g$-twisted $V$-module is admissible (cf. \cite{DLM1}).

The following result is a generalization of Proposition 4.5.7 of \cite{LL} (also see \cite{DLM3}, \cite{DLM4}, \cite{L3}, \cite{DJ}) and will be used later.

\begin{lem}\label{l5.2} Assume that $\mathbb{D}=\F.$ Let $W$ be a weak $g$-twisted $V$-module and $w\in W.$ Then for any $u\in V^r, v\in V^s$ and
$m\in \frac{r}{T}+\Z, n\in \frac{s}{T}+\Z$ there exist $m_1,...,m_k\in\Z,$ $n_1,...,n_k\in \frac{r+s}{T}+\Z$ and
$c_1,...,c_k\in \mathbb{D}$ such that
$$u_mv_nw=\sum_{i=1}^kc_i(u_{m_i}v)_{n_i}w.$$
In particular, the weak  $g$-twisted $V$-submodule of $W$ generated by $w$ is spanned by $a_nw$ for $a\in V$ and $n\in\frac{1}{T}\Z.$
\end{lem}

The proof is similar to that of Proposition 4.5.7 of \cite{LL} by using \ref{ea2.15} and noting that $T$ is invertible in the field $\F.$

We now define forms. Let $R$ be an  integral subdomain of $\mathbb{D}$ and $V$ a vertex operator superalgebra over $\mathbb{D}.$
\begin{definition} Assume that $V$ is a free $R$-module. Then $V$ is a vertex operator superalgebra over $R.$
A vertex operator superalgebra $R$-form $I$ of $V$ is a sub vertex operator superalgebra of $V$ such that
$I_n=I\cap V_n$  and $V_n=\mathbb{D}\otimes_RI_n$ for all $n\in \ha \Z.$
\end{definition}

Unlike \cite{DG} we do not assume that there is a nondegenerate symmetric bilinear form on $V.$ Also note that $I\cap V_{\bar{0}}$ is a vertex operator algebra $R$-form of $V_{\bar{0}}.$

\begin{definition} Let $V$ be a vertex operator superalgebra over $\mathbb{D}$ and $I$ a vertex operator superalgebra $R$-form of $V.$ Assume that $M=\oplus_{n\in\frac{1}{T}\Z_+}M(n)$ is an admissible $g$-twisted $V$-module
where $g$ is an automorphism of $V$ of order $T,$  $g$ preserves $I,$ $M$ is also an admissible $g$-twisted $I$-module. An $R$-form $N=\oplus_{n\in\frac{1}{T}\Z_+}N(n)$ of $M$ over $I$ is a $g$-twisted admissible $I$-submodule of $M$
such that $N(n)=N\cap M(n)$ and $M(n)=\mathbb{D}\otimes_RN(n)$ for all $n\in  \frac{1}{T}\Z_+.$
\end{definition}

We remark that if $I$ is a vertex operator superalgebra $R$-form of $V$, then $I\cap V_{\bar{1}}$ is an $R$-form of  $V_{\bar{1}}$ over $I\cap V_{\bar{0}}.$

\section{Vertex operator algebras associated to the Virasoro algebra}
\setcounter{equation}{0}
In this section we study the Virasoro algebra over an integral domain $\mathbb{D}$ and discuss the highest weight modules and related vertex operator algebras. In the rest of the paper we assume that $\F$ is a field with $\ch \F\ne 2.$

Recall that the original Virasoro algebra $Vir$ is an infinite dimensional Lie algebra over $\C$ with a basis
$L_n$ for $n\in\Z$ and $C.$  Consider the integral domain $\mathbb{D}$ such that $1/2\in \mathbb{D}.$ Set
$$Vir_\mathbb{D}=\oplus_{n\in \Z}\mathbb{D}L_n\oplus \mathbb{D}C$$
subject to the relation
$$[L_m,L_n]=(m-n)L_{m+n}+\frac{m^3-m}{12}\delta_{m+n,0}C, \  [Vir_\mathbb{D},C]=0$$
for $m,n\in \Z.$ Note that $m^3-m$ is divisible by $3$ for any $m\in \Z,$ the commutators make sense.
In particular we can define the Virasoro algebra $Vir_\F$ for any field $\F$ with $\ch \F\ne 2.$
In the rest of paper we always use $\F$ for such a field.  It is also clear that $Vir_\F=\F\otimes_RVir_R$ where $R=\Z[\ha].$

 As in the complex case, for any $c,h\in \mathbb{D}$ we set
$$V(c,h)_\mathbb{D}=U(Vir_\mathbb{D})\otimes_{U(Vir_\mathbb{D}^{\geq0})}\mathbb{D}$$
where $Vir_\mathbb{D}^{\geq0}$ is the subalgebra generated by $L_n$ for $n\geq 0$ and $C,$ and $\mathbb{D}$ is a
 $Vir_\mathbb{D}^{\geq0}$-module such that $L_n1=0$ for $n>0,$ $L_01=h$ and $C1=c.$ Then
  $$V(c,h)_\mathbb{D}=\oplus_{n\geq 0}V(c,h)_\mathbb{D}(n)$$
  is $\Z$-graded where $V(c,h)_\mathbb{D}(n)$ has a basis
  $$\{L_{-n_1}\cdots L_{-n_k}v_{c,h}|n_1\geq \cdots \geq n_k\geq 1, \sum_{i}n_i=n\}$$
  where $v_{c,h}=1\otimes 1.$ The $V(c,h)_\mathbb{D}$ is again is called the Verma module.
  It is easy to see that $L_nV(c,h)_\mathbb{D}(m)\subset V(c,h)_\mathbb{D}(m-n)$ for all $m,n\in\Z.$

 \begin{lem}\label{newl2.1} The $V(c,h)_\F$ has a unique maximal graded  submodule $W(c,h)_\F$ such that $L(c,h)_\F=V(c,h)_\F/W(c,h)_\F$ is an irreducible highest weight  $Vir_\F$-module. Moreover, if $\ch\F=0$
 then $W(c,h)_\F$ is the unique maximal submodule of $V(c,h)_\F.$
  \end{lem}

  \pf This result is easy if $\ch\F=0$  as each $V(c,h)_\F(n)$ is an eigenspace for $L_0$ with eigenvalue
  $h+n.$ If $\ch\F$ is finite , $L_0$ may have finitely many eigenvalues and
  the old proof does not work. This should explain why the statement of the lemma requires the maximal graded submodule.  We will show how to use the gradation to establish the result.

  Let $W(c,h)_\F=\oplus_{n\geq 1}W(c,h)_\F(n)$ is the unique maximal graded submodule of $V(c,h)_\F$ where
  $W(c,h)_\F(n)=W(c,h)_\F\cap V(c,h)_\F(n)$ for all $n.$ It is good enough to prove that $W(c,h)_\F$
is  a maximal submodule. Assume that
  $M$ is a submodule of $V(c,h)_\F$ strictly containing $W(c,h)_\F.$ We need to prove that $M=V(c,h)_\F.$ If not, $M$ contains a nonhomogeneous vector  $w=w^{n_1}+\cdots+ w^{n_k}$ which does not lie in $W(c,h)_\F,$ where $0\ne w^{n_i}\in V(c,h)_\F(n_i)\setminus W(c,h)_\F$ and $n_1>\cdots >n_k.$
  Now we consider the submodule $W$ of $V(c,h)_\F$ generated by $w^{n_1}.$ Since $w^{n_1}$ is a homogeneous vector we see that $W=\oplus_{n\geq 0}W(n)$ where $W(n)=W\cap V(c,h)_\F(n).$ If $W(0)=0$ then $W$ is contained in
  $W(c,h)_\F$ and $w^{n_1}\in W(c,h)_\F,$ a contradiction. This implies that $W(0)=V(c,h)_\F(0)\ne 0.$
  As a result there exists an element $a$  in $U(Vir_\F)$ which is a linear combination of
  $L(s_1)\cdots L(s_q)$ with $s_1\geq \cdots s_q\geq 1$ satisfying $\sum_i{s_i}=n_1$ such
  that $aw^{n_1}=v_{c,h}.$ It is clear that $aw^{n_j}=0$ for $j>1.$ So $aw=v_{c,h}\in M$ and $M=V(c,h)_\F,$ as expected.
  \qed

We now assume that $\ch\F=p$ is a prime and discuss an arbitrary maximal submodule $W$ of $V(c,h)_\F.$
Note that
$$V(c,h)_\F=\oplus_{r=0}^{p-1}V(c,h)^r_{\F}$$
where
$$V(c,h)_\F^r=\{v\in V(c,h)_\F|L_0v=(h+r)v\}=\oplus_{n\geq 0}V(c,h)_\F(r+pn).$$
 Then
$$W=\oplus_{r=0}^{p-1}W^r$$
where $W^r=W\cap V(c,h)_\F^r$ for all $r$ using the fact that $L_0$ has different enigenvalues on $V(c,h)_\F^r$
for different $r.$ As a result, the irreducible quotient module $V(c,h)_\F/W$ has the decomposition
$$V(c,h)_\F/W=\oplus_{r=0}^{p-1}V(c,h)_\F^r/W^r.$$
Clearly, if $W$ is different from $W(c,h)_\F$ then $W^0$ must contain a nonhomogeneous vector $w=\sum_{n\geq 0}w^{pn}$ with $w^{np}\in V(c,h)_\F(np)$ and $w^0\ne 0.$ It is definitely an interesting  problem to determine maximal submodules of $V(c,h)_\F.$ For any such $W,$ $V(c,h)_\F/W$ is not a $\Z$-graded module anymore as $W$ is not $\Z$-graded.

In the rest of paper, we will drop $\F$ from $Vir_\F$, $V(c,h)_\F,$ $\bar V(c,0)_\F$ and $L(c,h)_\F$ if
$\F=\C.$

We now turn our attention to $V(c,0)_\F.$ Set $\bar V(c,0)_\F=V(c,0)_\F/U(Vir_\F)L_{-1}v_{c,0}$ and denote the image of $v_{c,0}$ in $\bar V(c,0)_\F$ by $\1.$ Also let $\omega=L_{-2}\1\in \bar{V}(c,0)_\F.$ The following theorem
is an analogue of the same result in the complex case in the current setting with the same proof.

 \begin{thm}\label{t3.1} Let $c\ne 0.$ Then

(1) $(\bar{V}(c,0)_\F, Y, \1,\omega)$ is a vertex operator algebra generated by $\omega$ such that $Y(\omega, z)=
\sum_{n\in\Z}L_nz^{-n-2}.$

(2) For any $h,$ $V(c,h)_\F$ and its quotient modules for the Virasoro algebra are the $\bar{V}(c,0)_\F$-modules with $Y(\omega,z)=
\sum_{n\in\Z}L_nz^{-n-2}.$

(3) $W(c,0)_\F/U(Vir_\F)L_{-1}v_{c,0}$ is an ideal of $\bar{V}(c,0)_\F$ and $L(c,0)_\F$ is a simple vertex operator algebra which is a quotient vertex operator algebra of $\bar{V}(c,0)_\F$ modulo the ideal $W(c,0)_\F/U(Vir_\F)L_{-1}v_{c,0}.$

(4) Let $M$ be a restricted module for the Virasoro algebra $Vir_\F$ of central charge $c,$
that is, for any $w\in M$ $L_nw=0$ if $n$ is sufficiently large. Then $M$ is a weak  $\bar{V}(c,0)_\F$-module
with $Y(\omega, z)=\sum_{n\in\Z}L_nz^{-n-2}.$ Moreover, $M$ is an irreducible $Vir_\F$-module if and only if
$M$ is an irreducible $\bar{V}(c,0)_\F$-module.
\end{thm}

We next discuss the connection of $V(c,h)_\F, \bar V(c,0)_\F, L(c,h)_\F$ with the corresponding objects $V(c,h), \bar V(c,0), L(c,h)$ in the complex case via the  $R$-forms where $R=\Z[\ha]$ and  $c,h\in R.$ We also define $V(c,h)_R, \bar V(c,0)_R$ in the same way. Let $L(c,h)_R$ be the $R$-span of
$$\{L_{-n_1}\cdots L_{-n_k}(v_{c,h}+W(c,h))|n_1\geq \cdots \geq n_k\geq 1\}$$
in $L(c,h).$ Then it is clear $W_\F=\F\otimes_RW_R$ for $W=V(c,h), \bar V(c,0).$ It is natural to ask if $L(c,h)_\F$ is isomorphic to $\F\otimes_RL(c,h)_R?$ It is well known that $L(c,0)=\bar V(c,0)$ if $c>1$ \cite{KR}. But it is not true in the current situation.

\section{Modular vertex operator algebra $L(\ha,0)_\F$ }
\setcounter{equation}{0}

We investigate the modular vertex operator algebra $L(\ha,0)_\F$ in this section. In particular, we show that
$L(\ha,0)_\F$ is isomorphic to $\F\otimes_RL(\ha,0)_R$ by constructing $\F\otimes_RL(\ha,0)_R$ explicitly.

We first recall the construction of $L(\ha,0)$ and its irreducible modules $L(\ha,h)$ for $h=0,\ha,\se.$

Again let $\mathbb{D}$ be an integral domain such that $\ha\in \mathbb{D}.$  Let $H_\mathbb{D}=\mathbb{D} a$ be one dimensional free $\mathbb{D}$-module
with a nondegenerate bilinear form $(,)$ determined by $(a,a)=1.$ For short we also set $H=H_\C.$
For $Z=\Z$ or $\Z+\ha$ we let $A(H_\mathbb{D},Z)$ be an associative algebra over $\mathbb{D}$ generated by $a(n)$ for $n\in Z$ subject to the relation $a(m)a(n)+a(n)a(m)=\delta_{m+n,0}$ for $m,n\in Z.$ Also let $A(H_\mathbb{D},Z)_+$ be the subalgebra of $A(H_\mathbb{D},Z)$ genrated by $a(n)$ for $n>0$ and  $\mathbb{D}$ is  an $A(H_\mathbb{D},Z)_+$ -module such that $a(n)1=0$ for $n>0.$ Consdier the induced module $V(H_\mathbb{D},Z)=A(H_\mathbb{D},Z)\otimes_{A(H_\mathbb{D},Z)_+}\mathbb{D}.$ It is well known that $V(H_\F,\Z+\ha)$ is a simple $A(H_\F,\Z+\ha)$-module and $V(H_\F,\Z)$ is a direct sum of two simple $A(H_\F,\Z)$-modules.

Note that $V(H_\mathbb{D},Z)=\wedge[a(-n)|n\geq 0, n\in Z]$ as  $\mathbb{D}$-modules and the action of $A(H_\mathbb{D},Z)$ is as follows:
$a(n)$  acts as $\frac{\partial}{\partial a(-n)}$ for $n>0$ and as the multiplication by $a(n)$ for $n\leq 0.$
The $V(H_\mathbb{D},Z)$ decomposes into $V(H_\mathbb{D},Z)_{\bar 0}\oplus V(H_\mathbb{D},Z)_{\bar 1}$ where $V(H_\mathbb{D},Z)_{\bar s}$ is the submodule spanned by the monomials whose length is congruent to $s$ modulo $2.$
Set $\omega=\frac{1}{2}a(-\frac{3}{2})a(-\frac{1}{2})\in V(H_\mathbb{D},\Z+\ha).$
Let $\tau$ be the canonical automorphism of  $V(H_\mathbb{D}, \Z+\ha)$ so that $\tau|_{V(H_\mathbb{D},\Z)_{\bar s}}=(-1)^s.$
The following result is well known (cf. \cite{FFR},
 \cite{KW} and \cite{L1}, \cite{L2}, \cite{KR}).
 \begin{prop}\label{p4.1} (1) $(V(H_\mathbb{D}, \Z+\ha), Y, 1,\omega)$ is a vertex operator superalgebra generated by $a(-\ha)$ with  $$Y(a(-1/2),z)=\sum_{n\in\Z+\ha}a(n)z^{-n-\ha}$$
 and $V(H_\F, \Z+\ha)$ is a simple vertex operator superalgebra.

(2)  $(V(H_\mathbb{D}, \Z), Y)$ is a $\tau$-twisted  admissible $V(H_\mathbb{D}, \Z+\ha)$-module with $$Y(a(-1/2),z)=\sum_{n\in\Z}a(n)z^{-n-\ha}.$$

(3) $V(H,\Z+\ha)_{\bar 0}$ is isomorphic to $L(\ha,0),$ $V(H ,\Z+\ha)_{\bar 1}$ is isomorphic to $L(\ha,\ha),$
and $V(H,\Z)_{\bar s}$ is isomorphic to $L(\ha,\se)$ for $s=0,1$ as modules for the Virasoro algebra.

(4) $V(H_R,\Z+\ha)$ is an $R$-form of $V(H,\Z+\ha)$ and $V(H_R,\Z)$ is an $R$-form of $V(H,\Z)$ over $V(H,\Z+\ha)$
where $R=\Z[\ha].$
Moreover, $V(H_\F,Z)=\F\otimes_RV(H_R,Z)$ for any field $\F$ with $\ch\F\ne 2.$

(5) $V(H_R,\Z+\ha)_{\bar 0}$ is an $R$-form of  $V(H,\Z+\ha)_{\bar 0}$ and $V(H_R,Z)_{\bar s}$ is an $R$-form
of $V(H,Z)_{\bar s}$ over $V(H_R,\Z+\ha)_{\bar 0}$ for $Z=\Z+\ha,\Z$ and $s=0,1.$
\end{prop}

Here we give the  explicit expression for the $Y(\omega,z)=\sum_{n\in\Z}L(n)z^{-n-2}.$ For this purpose we define a normal ordering
$$:a(m)a(n):=\left\{\begin{array}{ll} a(m)a(n) & {\rm if}\  m\leq n\\
-a(n)a(m)  & {\rm if}\  m> n
\end{array}\right.$$
It is easy to check that
\begin{equation}\label{v1}
L(n)=\ha\sum_{j\in Z}j:a(-j)a(n+j):
\end{equation}
if $Z=\Z+\ha$ or $n\ne 0,$ and
\begin{equation}\label{v2}
L(0)=\frac{1}{16}+\ha\sum_{j\in \Z}j:a(-j)a(j):.
\end{equation}

The following lemma will be used later.
\begin{lem}\label{l4.2} For any field $\F,$ $V(H_\F,\Z+\ha)_{\bar 0}$ is a simple vertex operator algebra and $V(H_\F,Z)_{\bar s}$ for all $s$ are irreducible admissible $V(H_\F,\Z+\ha)_{\bar 0}$-modules.
Moreover, $V(H_\F,\Z)_{\bar 0}$ and $V(H_\F,\Z)_{\bar 1}$ are isomorphic.
\end{lem}

\pf First we show that $V(H_\F,Z)_{\bar s}$ is an admissible $V(H_\F,\Z+\ha)_{\bar 0}$-module. We set
$$V(H_\F,\Z+\ha)_{\bar s}(n)=\<a(-n_1)\cdots a(-n_k)\in V(H_\F,\Z+\ha)|n_1>\cdots >n_k>0,\sum_{i}n_i=n+\ha s\>$$
and
$$V(H_\F,\Z)_{\bar s}(n)=\<a(-n_1)\cdots a(-n_k)\in V(H_\F,\Z)_{\bar s} |n_1>\cdots >n_k>0,\sum_{i}n_i=n\>$$
for $s=0,1.$
Then
$$V(H_\F,Z)_{\bar s}=\oplus_{n\geq 0}V(H_\F,Z)_{\bar s}(n)$$
and
$$u_nV(H_\F,Z)_{\bar s}(m)\subset V(H_\F,Z)_{\bar s}(m+t-n-1)$$
for $u\in (V(H_\F,\Z+\ha)_{\bar 0})_t=V(H_\F,\Z+\ha)_{\bar 0}(t)$ for all $m,n,t.$

We next to prove each $V(H_\F,Z)_{\bar s}$  is an irreducible admissible $V(H_\F,\Z+\ha)_{\bar 0}$-module.
Note that $V(H_\F,Z)_{\bar s}$ is irreducible under the operators $a(m)a(n)$ for $m,n\in Z.$ Also note that
$a(m)=a(-\ha)_{m-\ha}.$ Let $w\in V(H_\F,Z)_{\bar s}.$ By Lemma \ref{l5.2} and Proposition \ref{p4.1} there exist $m_1,...,m_k\in\ha+\Z,$ $n_1,...,n_k\in \Z$ and
$c_1,...,c_k\in \F$ such that
$$a(m)a(n)w=\sum_{i=1}^kc_i(a(m_i)a(-1/2))_{n_i}w.$$
This shows that $V(H_\F,Z)_{\bar s}$ is an irreducible $V(H_\F,\Z+\ha)_{\bar 0}$-module. In particular,
$V(H_\F,\Z+\ha)_{\bar 0}$ is a simple vertex operator algebra.

It remains to show that  $V(H_\F,\Z)_{\bar 0}$ and $V(H_\F,\Z)_{\bar 1}$ are isomorphic. Define
a linear isomorphism $\sigma: V(H_\F,\Z)_{\bar 0}\to V(H_\F,\Z)_{\bar 1}$ such that
$a(-n_1)\cdots a(-n_k)$ for $n_1>\cdots>n_k\geq 0$ is sent to $a(-n_1)\cdots a(-n_k)a(0).$
From the discussion above, it is enough to verify that $a(s)a(t)\sigma=\sigma(s)a(t)$ for any $s,t\in\Z$ with $s> t.$ This is clear if both $s,t$ are different from $0.$ We now assume $t=0.$ Then
$$a(s)a(0)\sigma a(-n_1)\cdots a(-n_k)=\frac{\partial}{\partial a(-s)}a(0)a(-n_1)\cdots a(-n_k)a(0)$$
$$\sigma a(s)a(0)a(-n_1)\cdots a(-n_k)=\frac{\partial}{\partial a(-s)}a(0)a(-n_1)\cdots a(-n_k)a(0)$$
are equal for any $n_1>\cdots n_k\geq 0.$ The proof for $s=0$ is similar. \qed

We are now present our first main result of the paper which is a generalization of the same result from the complex field to any field $\F$ whose characteristic is not $2.$

\begin{thm}\label{tm1} If $\ch \F\ne 7,$ then $V(H_\F,\ha+\Z)_{\bar 0}$ is isomorphic $L(\ha,0)_\F,$
$V(H_\F,\ha+\Z)_{\bar 1}$ is isomorphic $L(\ha,\ha)_\F,$ and $V(H_\F,\Z)_{\bar s}$ is isomorphic $L(\ha,\se)_\F$
for $s=0,1$ as modules for the Virasoro algebra $Vir_\F.$
\end{thm}

\pf The results in the case that $\F=\C$ were proved using the classification of unitary highest weight modules for the Virasoro algebra with $c=\ha$ (cf. \cite{KR}). But this method does not work in the current situation.
We need to establish directly that $V(H_\F,Z)_{\bar s}$ is an irreducible highest weight module for the Virasoro algebra $Vir_\F.$ Of course, the proof present here also works for $\F=\C.$

Note that $\1\in V(H_\F,\Z+\ha)_{\bar 0},$ $a(-1/2)\in V(H_\F,\Z+\ha)_{\bar 1},$
$\1\in V(H_\F,\Z)_{\bar 0},$ and $a(0)\in V(H_\F,\Z)_{\bar 1}$  are highest weight vectors for $Vir_\F$ with highest weight $0,\ha, \se,\se.$ This is true for $\F=\C$ and the same argument works for any $\F$ with $\ch\F\ne 2.$ So if we can prove that $V(H_\F,Z)_{\bar s}$ is an irreducible $Vir_\F$-module, then the results follow.

We first show that $V(H_\F,\Z+\ha)_{\bar 0}$ is an irreducible module for $Vir_\F.$
Let $W$ be the $Vir_\F$-submodule of $V(H_\F,\Z+\ha)_{\bar 0}$ generated by $\1.$ From the proof of Lemma \ref{l4.2} we only need to argue that $a(-s-\ha)a(-\ha)\in W$ for all $0<s\in\Z.$ We prove by induction on
$s$ that $a(-m-\ha)a(-n-\ha)$ lie in $W$ for $m,n\geq0$ with $m+n=s.$  We need the following commutator relation
$$[L(p),a(q)]=-(q+\frac{p}{2})a(p+q)$$
which is a consequence of (\ref{g2.17}) for $p\in\Z$ and $q\in\frac{1}{2}+\Z.$

If $s=1$ then $a(-3/2)a(-1/2)=2L(-2)\1=2\omega\in W.$ If $s=2$ then $a(-5/2)a(-1/2)=L(-1)\omega\in W.$ If $s=3$
consider the linear system
$$L(-1)a(-\frac{5}{2})a(-\ha)=3a(-\frac{7}{2})a(-\ha)+a(-\frac{5}{2})a(-\frac{3}{2})$$
$$L(-2)a(-\frac{3}{2})a(-\ha)=(2+\ha)a(-\frac{7}{2})a(-\ha)-\frac{3}{2}a(-\frac{5}{2})a(-\frac{3}{2}).$$
Since the coefficient matrix
$$\left(\begin{array}{cc} 3 & 1\\ 2+\ha & -\frac{3}{2}\end{array}\right)$$
has determinant
$-7$ we see that both $a(-\frac{7}{2})a(-\ha)$ and $a(-\frac{5}{2})a(-\frac{3}{2})$ are in $W$ as $\ch\F\ne 7.$

We now assume that the result holds for $s-1$ and $s$.
We have linear system
$$L(-1)a(-s+i-\ha)a(-i-\ha)=(s-i+1)a(-s+i-\frac{3}{2})a(-i-\ha)+(i+1)a(-s+i-\ha)a(-i-\frac{3}{2})$$
%$$L(-1)a(-s-\ha)a(-\ha)=(s+1)a(-s-\frac{3}{2})a(-\ha)+a(-s-\ha)a(-\frac{3}{2})$$
%$$L(-1)a(-s+\frac{1}{2})a(-\frac{3}{2})=sa(-s-\frac{1}{2})a(-\frac{3}{2})+2a(-s+\frac{1}{2})
%a(-\frac{5}{2})$$
$$L(-2)a(-s+\ha)a(-\ha)=(s+\ha)a(-s-\frac{3}{2})a(-\ha)+\frac{3}{2}a(-s+\ha)a(-\frac{5}{2})$$
$$L(-2)a(-s+\frac{3}{2})a(-\frac{3}{2})=(s-\ha)a(-s-\frac{1}{2})a(-\frac{3}{2})+\frac{5}{2}a(-s+\frac{3}{2})
a(-\frac{7}{2})$$
%$$L(-3)a(-s+\frac{3}{2})a(-\ha)=sa(-s-\frac{3}{2})a(-\ha)+2a(-s+\frac{3}{2})
%a(-\frac{7}{2})$$
for $i=0,...,[s/2].$
By induction assumption, the vectors in the left hand side of the linear system lie in $W.$ %The coefficient matrix is
%$$\left(\begin{array}{cccc} s+1 & 1 & 0 &0\\ 0& s& 2&0\\ s+\ha & 0& \frac{3}{2}&0\\ %0&s-\ha&0&\frac{5}{2}\end{array}\right)$$
Solving the linear system above shows that $a(-m-\ha)a(-n-\ha)$ lies in $W$  for $m+n=s+1.$ The proof is complete.

We now prove the others. By Lemma \ref{l4.2}, each $V(H_\F,Z)_{\bar s}$ is an irreducible $V(H_\F,\Z+\ha)_{\bar 0},$ or an irreducible $L(\ha,0)_\F$-module. Theorem \ref{t3.1} tells us that $V(H_\F,Z)_{\bar s}$ is an
irreducible $Vir_\F$-module. \qed

The following corollary gives the degrees of the singular vectors of $V(\ha,h)_\F$ for $h=0,\ha,\se.$
\begin{cor}\label{c1} Assume that $\ch \F\ne 7.$ Then the maximal graded submodule $W(\ha,h)_\F$ of $V(\ha,h)_\F$ is generated by two singular vectors which have degrees $1$ and $6$ if $h=0,$
$2$ and $3$ if $h=\ha,$ and $2$ and $4$ if $h=\se.$
\end{cor}

\pf We first recall the expressions of the singular vectors
of $V(\ha,h)$ for $h=0,\ha,\se$ from \cite{FF} (also see \cite{BFIZ}).
The $V(\ha,0)$ has two singular vectors
$$L(-1)v_{\ha,0},\ (64L(-2)^3+93L(-3)^2-264L(-4)L(-2)-108L(-6))v_{\ha,0}$$
of degrees $1$ and $6,$
$V(\ha,\ha)$ has two singular vectors
$$(4L(-2)-3L(-1)^2)v_{\ha,\ha}, (L(-1)^3-3L(-2)L(-1)+\frac{3}{4}L(-3))v_{\ha,\ha}$$
of degrees $2$ and $3,$ $V(\ha,\se)$ has two singular vectors
of degrees $2$ and $4$
$$(3L(-2)-4L(-1)^2)v_{\ha,\se},$$
$$(L(-1)^4+2uL(-2)L(-1)^2+(u^2-4)L(-2)^2+(2u+6)L(-3)L(-1)+(3u+6)L(-4))v_{\ha,\se}$$
where $u=-\frac{25}{12}$ and $L(n)$ is identified with $L_n$ for $n\in\Z.$

Since all coefficients of these singular vectors are in $\Z[\ha]$ we
see that they are also the singular vectors $V(\ha,h)_\F$ for any field $\F.$
By Theorem \ref{tm1}, $L(\ha,h)_\F$ have the same graded dimensions for all fields $\F$ with $\ch\F\ne 7.$
Since the result is true for $\F=\C$ from \cite{FF}, the results also hold for any field $\F.$ \qed

Recall $R=\Z[\ha]$. It follows from Proposition \ref{p4.1} that $L(\ha,0)_R=V(H_R,\Z+\ha)_{\bar 0}$ is an $R$-form of $L(\ha,0),$ $L(\ha,\ha)_R=V(H_R,\Z+\ha)_{\bar 1}$ is an $R$-form
of $L(\ha,\ha)$ over $L(\ha,0)_R$ and $L(\ha,\se)_R=V(H_R,\Z)_{\bar s}$ is an $R$-form of $L(\ha,\se)$ over $L(\ha,0)_R.$
Theorem \ref{tm1}  gives the following corollary.
\begin{cor}\label{c2} If $\ch\F\ne 7$ then $L(\ha,h)_\F=\F\otimes_RL(\ha,h)_R$ for $h=0,\ha,\se.$
\end{cor}

We next classify the irreducible admissible modules for $L(\ha,0)_{\F}.$ It is well known that $L(\ha,0)$ has only three irreducible admissible modules $L(\ha,h)$ with $h=0,\ha,\se$ (see \cite{DMZ} and \cite{W}). Although the proof in the complex case uses the Zhu's algebra $A(V)$ \cite{Z}, the key idea is to find the image of the singular vectors in $A(V)$. But there is no Zhu's algebra theory for vertex operator algebra over arbitrary field $\F$ currently (Zhu's $A(V)$-theory is studied in a forthcoming paper), we will get the classification result directly using the singular vectors.

\begin{thm}\label{tm2} Let $M$ be a highest weight module for $Vir_\F$ with central charge $\ha$ and highest weight $h\in\F$ such that $M$ is also a module for $L(\ha,0)_\F$ and $M(0)\ne 0.$
Then $h=0,\ha,\se$ if $\ch \F\ne 7$ and
$h=0, \ha$ if $\ch\F=7.$ In particular, the vertex operator algebra $L(\ha,0)_\F$ has  three irreducible admissible modules $L(\ha,h)_\F$  with $h=0,\ha,\se$   if $\ch \F\ne 7$ and has two possible irreducible admissible modules are $L(\ha,h)_\F$ with $h=0,\ha$ if  $\ch \F= 7$ up to isomorphism.
\end{thm}

\pf  It is well known that $\bar V(\ha,0)$ has a unique maximal submodule generated by a highest weight vector ${\bf s}$ of weight $6$ when $\F=\C.$ The vector ${\bf s}$ is given by $${\bf s}=(64L(-2)^3+93L(-3)^2-264L(-4)L(-2)-108L(-6))\1$$
\cite{FF}. Regarding the coefficients appearing in the expression of ${\bf s}$  as elements in $\F$ we see that ${\bf s}$ is also a singular  vector of ${\bar V}(\ha,0)_{\F}.$
This implies that ${\bf s}$ is  the zero vector in
$L(\ha,0)_{\F}$  and $Y({\bf s},z)=0$ on any weak $L(\ha,0)_\F$-module. In particular, ${\bf s}_n=0$ on $M$ for all $n\in \Z.$

Let $V$ be a vertex operator algebra over $\F$ and $W$ a weak $V$-module. Then we have the following relation
on $W$ for any $u\in V$ by using the Jacobi identity
\begin{equation}\label{e4.1}
Y(u_{-1}v,z)=Y(u,z)^-Y(v,z)+Y(v,z)Y(u,z)^+
\end{equation}
where
$$Y(u,z)^+=\sum_{n\geq 0}u_nz^{-n-1}, Y(u,z)^-=\sum_{n<0}u_nz^{-n-1}$$
(see \cite{DL}). As a result we see that
\begin{eqnarray*}
& &\ \ \ Y(u_{-1}u_{-1}u,z)\\
& &=Y(u,z)Y(u,z)^+Y(u,z)^++2Y(u,z)^-Y(u,z)Y(u,z)^++ Y(u,z)^-Y(u,z)^-Y(u,z).
\end{eqnarray*}

Now let $0\ne v\in M(0).$ Using (\ref{e4.1}) we have
\begin{eqnarray*}
& &(L(-2)^3\1)_{5}v=(L(1)L(-1)L(0) +L(0)L(0)L(0)\\
& &\ \ \ +L(2)L(-1)L(-1)+L(1)L(0)L(-1)+L(0)L(1)L(-1))v\\
& &\ \ =(h^3+6h^2+8h)v,
\end{eqnarray*}
\begin{eqnarray*}
(L(-3)^2\1)_{5}v=(4L(0)L(0)+3L(1)L(-1))v=(4h^2+6h)v,
\end{eqnarray*}
\begin{eqnarray*}
(L(-4)L(-2)\1)_{5}v=(3L(0)L(0)+L(1)L(-1))v=(3h^2+2h)v,
\end{eqnarray*}
$$(L(-6)\1)_5v=5hv.$$
So we have
$$0={\bf s}_{5}v=64h(h-\ha)(h-\se)v.$$
That is, $h=0, \ha,\se.$ So the possibilities of  irreducible admissible $L(\ha,0)_\F$-modules
are  $L(\ha,h)_\F$  with $h=0,\ha,\se$ if $\ch\F\ne 7$ and with $h=0,\ha$ if $\ch\F=7.$ By Theorem \ref{tm1}, those modules are indeed $L(\ha,0)_{\F}$-modules if $\ch\F\ne 7.$
\qed

It is worthy to point out that if $\ch\F=7,$ we will find a degree 4 singular vector which will give exactly the condition $h(h-1/2)$ for the any irreducible highest weight module $L(\ha,h)_\F.$

We next deal with the rationality of $L(\ha,0)_\F.$ In the case that $\F=\C,$ this is essentially proved in \cite{FF} (also see \cite{DMZ}, \cite{W}). Establishing rationality for general field $\F$ is  much more difficult. For example, it is almost trivial to prove that any admissible $L(\ha,0)$-module which is a highest weight module for the Virasoro algebra $Vir$ is an irreducible as the difference of any two numbers in $\{0,\ha,\se\}$ is not an integer. But this is not the case anymore if the characteristic is finite. The following lemma tells the same result is still true.
\begin{lem}\label{l2} If $M$ is $L(\ha,0)_\F$-module which is also a highest weight module for $Vir_\F,$ then
$M$ is irreducible and necessarily isomorphic to $L(\ha,h)_\F$ for $h=0,\ha,\se.$
\end{lem}

\pf Let $h$ be the highest weight of $M.$ Then $M$ is a quotient of $V(\ha,h)_\F$ generated by a highest weight vector $v$ with $h=0,\ha,\se$ by Theorem \ref{tm2}.

We first consider the case $h=0.$ As in the proof of Theorem \ref{tm2} we have the following computation
\begin{eqnarray*}
& &(L(-2)^3\1)_{6}L(-2)v=(L(0)L(-1)L(2)+L(-1)L(0)L(2)+L(1)L(-1)L(1)\\
& &\ \ \ +L(0)L(0)L(1)+L(-1)L(1)L(1)+L(2)L(-1)L(0)+L(1)L(0)L(0)\\
& &\ \ \ +L(0)L(1)L(0)+L(-1)L(2)L(0)+L(3)L(-1)L(-1)+L(2)L(0)L(-1)\\
& &\ \ \ +L(1)L(1)L(-1)+L(0)L(2)L(-1)+L(-1)L(3)L(-1))L(-2)v\\
& &\ \ =(121+\frac{3}{4})L(-1)v,
\end{eqnarray*}
\begin{eqnarray*}
& &(L(-3)^2\1)_{6}L(-2)v=(4L(-1)L(2)+6L(0)L(1)+6L(1)L(0)\\
& &\ \ \ +4L(2)L(-1))L(-2)v\\
& &\ \ =92L(-1)v,
\end{eqnarray*}
\begin{eqnarray*}
& &(L(-4)L(-2)\1)_{6}L(-2)v=(10L(-1)L(2)+6L(0)L(1)+3L(1)L(0)\\
& &\ \ \ +L(2)L(-1))L(-2)v\\
& &\ \ =(47+\frac{3}{4})L(-1)v,
\end{eqnarray*}
$$(L(-6)\1)_6L(-2)v=45L(-1)v.$$
Thus
$$0={\bf s}_7L(-2)v=(64(121+\frac{3}{4})+93\cdot 92-264(47+\frac{3}{4})-108\cdot 45)L(-1)v$$
and $L(-1)v=0.$ Since $L(\ha,0)_\F$ is generated by $\omega,$ using (\ref{e4.1}) we see that
$u_nv=0$ for all $u\in L(\ha,0)_\F$ and $n\geq 0.$ That is, $v$ is a vacuum-like vector.
Using the same proof
of Proposition 4.7.7 of \cite{LL} shows that the map $f:L(\ha,0)_\F\to M$ by sending $u$ to $u_{-1}v$
is a $L(\ha,0)_\F$-module homomorphism. Since $M$ is a $L(\ha,0)_\F$-module generated by $v$ we see that
the map $f$ is an isomorphism. That is, $M$ is irreducible.

In fact, if the characteristic of $\F$ is different from $3,5$ we do not need the computations above.
Since  $L(-1)v$ is a highest weight module for the Virasoro algebra of weight $1$ which is
different from $\ha,\se,$ Theorem \ref{tm2} then forces $L(-1)v=0.$

We next deal with the case $h=\ha, \se.$ Recall the expressions of the singular vectors
of $V(\ha,\ha)$ and $V(\ha,\se)$ from the proof of Corollary \ref{c1}.  As before we have
\begin{eqnarray*}
& &(L(-2)^3\1)_{5}L(-2)v=(L(-1)L(-1)L(2)+3L(-2)L(0)L(2)+L(0)L(-1)L(1)\\
& &\ \ \ +L(-1)L(0)L(1)+3L(-2)L(1)L(1)+L(1)L(-1)L(0) +L(0)L(0)L(0)\\
& &\ \ \ +L(-1)L(1)L(0)+3L(-2)L(2)L(0)+ L(2)L(-1)L(-1)+L(1)L(0)L(-1)\\
& &\ \ \ +L(0)L(1)L(-1)+L(-1)L(2)L(-1)+3L(-2)L(3)L(-1) )L(-2)v\\
& &\ \ =p(h)L(-2)v+q(h)L(-1)^2v
\end{eqnarray*}
where
$$p(h)=(6h+7)(h+\frac{1}{4})+(h+2)(11h+24)+18h,$$
$$q(h)=2(h+\frac{1}{4})+18h+54,$$
\begin{eqnarray*}
& &(L(-3)^2\1)_{5}L(-2)v=(3L(-1)L(1)+4L(0)L(0)+3L(1)L(-1))L(-2)v\\
& &\ \ =((4(h+2)^2+6(h+2))L(-2)v+18L(-1)^2v,
\end{eqnarray*}
\begin{eqnarray*}
& &(L(-4)L(-2)\1)_{5}L(-2)v=(10L(-2)L(2)+6L(-1)L(1)+3L(0)L(0)\\
& &\ \ \ +L(1)L(-1))L(-2)v\\
& &\ \ =(10(\frac{1}{4}+h)+(h+2)(3h+10))L(-2)v+21L(-1)^2v,
\end{eqnarray*}
$$(L(-6)\1)_5L(-2)v=5(h+2)L(-2)v.$$
Again we have
$$0={\bf s}_{5}L(-2)v=f(h)L(-2)v+g(h)L(-1)^2v$$
where
$$f(h)=64p(h)+93(4(h+2)^2+6(h+2))-264(10(\frac{1}{4}+h)+(h+2)(3h+10))-108\cdot 5(h+2)$$
and
$$g(h)=64q(h)+93\cdot 18-264\cdot 21.$$
One can easily verify that ${\bf s}_{5}L(-2)v$ is a nonzero multiple of
$(4L(-2)-3L(-1)^2)v$ if $h=\ha$ and a nonzero multiple of $(3L(-2)-4L(-1)^2)v$ if $h=\se.$
Thus $(4L(-2)-3L(-1)^2)v=0$ for $h=\ha$ and $(4L(-2)-3L(-1)^2)v=0$ for $h=\se.$
Similarly, $(L(-1)^3-3L(-2)L(-1)+\frac{3}{4}L(-3))v=0$ for $h=\ha$ by using ${\bf s}_5L(-3)v=0$
and the singular vector of degree $4$ in $M$ is zero for $h=\se$ by using ${\bf s}_5L(-4)v=0.$

By Corollary \ref{c1}, $M$ is indeed an irreducible module for the Virasoro algebra $Vir_\F$ and for vertex operator algebra $L(\ha,0)_\F.$ The proof is complete.\qed

We remark that if $\ch\F>13$ or $\ch\F=0$ then $0, 6, \ha, \ha+2,\ha+3, \se, \se+2,\se+4$ are all different.
In these cases we immediately conclude that $M$ is irreducible by Corollary \ref{c1}. So the proof in Lemma \ref{l2} is  basically for $\ch\F=3,5,11,13$ although the proof works in general.
\begin{rem} Let $V$ be a vertex operator algebra and $W$ a weak $V$-module. It is known from \cite{LL} that if $v$ is a vacuum-like vector when $\ch\F=0$ if and only if $L(-1)v=0.$ But this result is not valid anymore if $\ch\F=p$ is a prime. For example, $L(-1)L(-np+1)\1=0$ for any positive integer $n$ but $L(-np+1)\1$ is not a vacuum-like vector
as $L(0)L(-np+1)\1=-L(-np+1)\1.$
\end{rem}

We are now in a position to establish the rationality of $L(\ha,0)_\F.$
\begin{thm}\label{tm3} If $\ch \F\ne 7,$ $L(\ha,0)_\F$ is rational.
\end{thm}

\pf Let $M=\oplus_{n\geq 0}M(n)$ be an admissible $L(\ha,0)_\F$-module such that $M(0)\ne 0.$ Let $W_n$ be the $L(\ha,0)_\F$-submodule generated by $M(s)$ for $s\leq n.$ Note that $W_n=\oplus_{s\geq 0}W_n(s)$ where
$W_n(s)=W_n\cap M(s).$ We now prove by induction on $n$ that
$W_n$ is completely reducible.

If $n=0,$ by Lemma \ref{l2} and the fact that $L(0)(L(0)-\ha)(L(0)-\se)v=0$ for any $v\in M(0)$ (see the proof
of Theorem \ref{tm2}), $W_0$ is completely reducible. We now assume that $W_s$ is completely reducible
for $s\leq n-1.$ If $W_{n-1}(n)=M(n)$ we have nothing to prove. So we assume that $W_{n-1}(n)\ne M(n).$
Consider the quotient admissible module
$$M/W_{n-1}=\oplus_{s\geq n}(M/W_{n-1})(s)=\oplus_{s\geq n}M(s)/W_{n-1}(s).$$
Let $v\in M(n)\setminus W_{n-1}(n).$ Clearly, $v+W_{n-1}$ is a highest weight vector in $M/W_{n-1}.$ Let $X$ be
a $L(\ha,0)_\F$-submodule of $M$ generated by $v.$ Note that $X=\oplus_{s\geq 0}X(s)$ where $X(s)=X\cap M(s)$
is finite dimensional. Then $Z=X\cap W_{n-1}$ is a submodule of $X$ and $X/Z$ is an irreducible
$L(\ha,0)_\F$-module. Note that $Z$ is completely reducible. We have an exact sequence
\begin{equation}\label{exact}
0\to Z\to X\to X/Z\to 0.
\end{equation}

Let $X'=\oplus_{s\geq 0}X(s)^*$ be the graded dual of $X.$ In the case $\ch\F=0,$ $X'$ is exactly the contragredient module of $X$ \cite{FHL}. The definition of contragredient module uses the operator
$e^{zL(1)}$ on vertex operator algebra. But the operator $e^{zL(1)}$ may not make sense if $\ch\F\ne 0.$ Nevertheless, we prove  that $X'$ is an admissible $L(\ha,0)_\F$-module without fully using the properties
of contragredient module.

Denote the natural non-degenerate paring $X'\times X\to \F$ by $(,).$ Clearly, $X'$ is a $Vir_\F$-module such that
$$(L(n)w', w)=(w',L(-n)w), (Cw',w)=(w',Cw)$$
for $n\in \Z, w\in X$ and $w'\in X'.$  Since $L(n)X(s)^*\subset X(s-n)^*$ for all $n,s$ we see that $X'$ is a restricted module for $Vir_\F$-module. By Theorem \ref{t3.1}, $X'$ is an admissible module for
$\bar V(\ha,0)_\F.$ To prove that $X'$ is an admissible module for
$L(\ha,0)_\F$ we need to verify that $Y({\bf s},z)=0$ on $X'$ by Corollary \ref{c1}.
A straightforward computation using (\ref{e4.1}) gives $({\bf s}_nw',w)=(w', {\bf s}_{10-n}w)=0$
for all $n\in\Z.$ One can also use the well known formula for the contragredient module in the complex case
$$(Y(u,z)w',w)=(-1)^{\deg u}z^{-2 {\deg u}}(w',Y(u,z^{-1})w)$$
in the current situation as long as $L(1)u=0.$ Since ${\bf s}$ is a highest weight vector in $\bar V(\ha,0)$
for the Virasoro algebra $Vir_\F,$ we immediately see that $( Y({\bf s},z)w',w)=0.$

Since each homogeneous subspace $X(s)$ of $X$ is finite dimensional, We have an exact sequence of admissible $L(\ha,0)_\F$-modules from (\ref{exact}):
$$0\to (X/Z)'\to X'\to Z'\to 0.$$
Note that $Z'$ is completely reducible $L(\ha,0)_\F$-module. We claim that $Z'$ is isomorphic to the  submodule $W$ of $X'$ generated by $Z(s)^*$ for $s\leq n-1.$ Assume that  $Z$ is generated by linearly independent
 highest weight vectors $w_i\in Z(s_i)$ of highest weights $h_i$ with $s_i<n$ for $i=1,...,k.$ Then both
 $Z$ and $Z'$  are isomorphic to $\oplus_{i=1}^kL(\ha,h_i)_\F$ as $L(\ha,0)_\F$-modules. Clearly, $W$ is also
 generated by a set of linearly independent highest weight vectors $w_i'\in Z(s_i)^*$ for $i=1,...,k.$ As a result, $W$ is isomorphic to $Z'$ and $Z'$ can be regraded as a submodule of $X'.$ That is,
 $X'$ is isomorphic to $(X/Z)'\oplus Z'.$ This implies that $X=X/Z\oplus Z$ is completely reducible.
 So the submodule of $M$ generated by any vector $v\in M(n)$ is completely reducible and $W_n$ is completely reducible. \qed

\section{The case $\ch\F=7$}
\setcounter{equation}{0}

 In this section we study $L(\ha,0)_\F$ when $\ch\F=7.$ According to Theorem \ref{tm2}, $L(\ha,0)_\F$ has two possible inequivalent irreducible admissible modules $L(\ha,0)_\F$ and $L(\ha,\ha)_\F.$   It is shown in this section that $L(\ha,0)_\F$ and $L(\ha,\ha)_\F$  are indeed $L(\ha,0)_\F$-modules. But we cannot determine the singular vectors for $V(\ha,h)_\F$ for $h=0,\ha$ and establish the rationality of $L(\ha,0)_\F.$

Recall the singular vector ${\bf s}$ from Section 4. A straightforward verification gives
 \begin{lem}\label{l6.1} The ${\bf u}=(L(-2)^2-2L(-4))\bf 1$ is a singular vector of $\bar V(\ha,0)_\F$ and  $$(L(-2)+L(-1)^2){\bf u}={\bf s}.$$
 \end{lem}

Let $M$ be a highest weight $Vir_\F$-module generated by a highest weight vector $w$ with highest weight $h.$ Assume that $M$  is also an $L(\ha,0)_\F$-module. It is easy to show that $0={\bf u}_3w=h(h-4)w.$ This result, of course, is consistent with Theorem \ref{tm2}.

Lemma \ref{l6.1} also tells us that $V(H_\F,\Z+\ha)_{\bar 0}$ is not equal to $L(\ha,0)_\F.$ To see this we
 recall that $V(H,\Z+\ha)_{\bar 0}\cong L(\ha,0)$ and  $V(\ha,0)$ does not have a singular vector of weight 4. On the other hand,  $V(H_\F,\Z+\ha)$ and $V(H,\Z+\ha)$ have the same graded dimensions. This implies that
 $V(H_\F,\Z+\ha)$ should have a highest weight vector of degree 4 corresponding to ${\bf u}.$ One can easily
 check that $a(-\ha)a(-\frac{7}{2})-3a(-\frac{3}{2})a(-\frac{5}{2})$ is a highest weight vector for the Virasoro algebra $Vir_\F.$

 \begin{thm}\label{tm4} The vertex operator algebra  $L(\ha,0)_\F$ has exactly two inequivalent irreducible admissible modules $L(\ha,0)_\F,$ $L(\ha,\ha)_\F.$
 \end{thm}

 \pf We prove the result by realizing $L(\ha,0)_\F$ as a subalgebra $V(H_\F,\Z+\ha)_{\bar 0}.$ We first show that
 $V(H_\F,Z)$ is a completely reducible $Vir_\F$-module.

 Define a nondegenerate symmetric bilinear form $(,)$ on $V(H_\F,Z)$ such that
 the monomials $a(-n_1)\cdots a(-n_k)$ form an orthonormal basis. Then the vectors
 from the different degrees are orthogonal and
 $$(a(n)u,v)=(u,a(-n)v)$$
 for all $n\in Z,$ and $u,v\in V(H_\F,Z)$ (cf. \cite{KR}).
 Recall from (\ref{v1}) and (\ref{v2}) the expression
 of $L(n)$ in terms of $a(s)a(t).$  Then
 \begin{equation}\label{v3}
 (L(n)u,v)=(u,L(-n)v).
 \end{equation}

 Notice that the restriction of the bilinear form to $V(H_\F,Z)_{\bar s}$ is nondegenerate in all cases. We now prove that $V(H_\F,\Z+\ha)_{\bar 0}$ is a completely reducible module for $Vir_\F$ and the proof of the other cases are identical.

Let $W_0$ be the $Vir_\F$-submodule of $V(H_\F,\Z+\ha)_{\bar 0}$ generated by $\1$ Then $W_0$ is a highest weight
 module for $Vir_\F$ with highest weight $0.$ Denote the maximal graded submodule of $W_0$ by $X_0.$
 Let $M_0$ be the maximal graded submodule of $V(H_\F,\Z+\ha)_{\bar 0}$ such that $M_0(0)=M_0\cap V(H_\F,\Z+\ha)_{\bar 0}(0)=0.$ Then $(\1, M_0)=0,$ $X_0\subset M_0$ and
 $V(H_\F,\Z+\ha)_{\bar 0}/M_0\cong L(\ha,0)_\F$ is irreducible. Moreover, $V(H_\F,\Z+\ha)_{\bar 0}=W_0+M_0.$
 Note that any element of $W_0$ is a linear combination of vectors of the form
 $$L(-n_1)\cdots L(-n_k)\1$$
 with  $n_1>\cdots n_k>0.$ Since
 $$(L(-n_1)\cdots L(-n_k)\1,M_0)=(\1,L(n_k)\cdots L(n_1)M_0)=0$$
 we see that $(W_0,M_0)=0.$ This implies that $(X_0, W_0+M_0)=(X_0,V(H_\F,\Z+\ha)_{\bar 0})=0$ and
 $X_0=0.$ As a result, $W_0=W_0/X_0\cong L(\ha,0)_\F$ is irreducible, the restriction of the bilinear form
 to $W_0$ is nondegenerate and $V(H_\F,\Z+\ha)_{\bar 0}=W_0\oplus M_0.$  In particular, $M_0=(W_0)^{\perp}$ and
the restriction of the bilinear form to $M_0$ is nondegenerate.

 We now let $W_1$ be the $Vir_\F$-submodule of $M_0$ generated by the $M_0(n_0)$ where $n_0$ is the smallest positive integer such that $M_0(n_0)\ne 0.$ Let $X_1$ be the maximal graded submodule of $W_1$ such that
 $X_1(n_0)=X_1\cap W_1(n_0)=X_1\cap M_0(n_0)=0$ and
 $M_1$ is the maximal graded submodule of $M_0$ such that $M_1(n_0)=M_1\cap M_0(n_0)=0.$ Then
 $X_1\subset M_1.$ Using the same proof we show that $X_1=0,$ $W_2$ is completely reducible $Vir_\F$-module
 and $M_1=(W_1)^{\perp}.$ Continuing this procedure asserts that $V(H_\F,\Z+\ha)_{\bar 0}$ is a completely reducible module for $Vir_\F.$

Thus, the $Vir_\F$-submodule of $V(H_\F,\Z+\ha)_{\bar 1}$ generated by $a(-\ha)$ is isomorphic to
$L(\ha,\ha)_\F$ and the $Vir_{\F}$-submodules of $V(H_\F,\Z)$ generated by ${\bf 1}$ and $a(0)$
are isomorphic to $L(\ha,\se)_\F=L(\ha,\ha)_\F.$ It is clear now that both $L(\ha,0)_\F$ and $L(\ha,\ha)_\F$
are $L(\ha,0)_\F$-modules.
\qed

From Theorem \ref{tm4} we see that $\F\otimes_R L(\ha,0)_R$ is a completely reducible module for $L(\ha,0)_\F.$ We have already known that  $\F\otimes_R L(\ha,0)_R$ is not a simple vertex operator algebra anymore. In fact, $\F\otimes_R L(\ha,0)_R$ has a highest weight vector $a(-\ha)a(-\frac{7}{2})-3a(-\frac{3}{2})a(-\frac{5}{2}).$ One can also verify that
  $$-a(-\frac{15}{2})+a(-\ha)a(-\frac{3}{2})a(-\frac{11}{2})
  +a(-\ha)a(-\frac{5}{2})a(-\frac{9}{2})+3a(-\frac{3}{2})a(-\frac{5}{2})a(-\frac{7}{2})$$
  is a highest weight vector of $\F\otimes_R L(\ha,\ha)_R.$ This implies that $\F\otimes_RL(\ha,\ha)_R$ is not irreducible module for $Vir_\F$ or $L(\ha,0)_\F.$ It is definitely an interesting question
  to find the highest weight vectors in $\F\otimes_R L(\ha,h)_R$ for $h=0,\ha,\se.$ This will eventually lead to the character formula for $L(\ha,h)_\F$ and  expressions of the singular vectors of $V(\ha,h)_\F$ for $\ch\F=7.$


\begin{thebibliography}{ABCDE}
\bibitem[BFIZ]{BFIZ} M. Bauer, ph. Di Francesco, C. Itzykson and J. B. Zuber, Covariant differential equations and singular vectors in Virasoro representations,
{\em Nucl. Phys. } {\bf B362} (1991), 515-562.

\bibitem[B]{B}R. E. Borcherds, Vertex algebras, Kac-Moody algebras, and the Monster,
{\em Proc. Natl. Acad. Sci. USA} {\bf 83} (1986), 3068-3071.

\bibitem[D]{D} C. Dong, Representations of the moonshine module vertex operator
algebra, {\em Contem. Math.}, {\bf 175} (1994), 27--36.

\bibitem[DG]{DG}C. Dong and R. L. Griess, Integral forms in vertex operator algebras, {\bf 365} {\em J. Algebra}  {\bf 365} (2012), 184-198.

\bibitem[DGH]{DGH}C. Dong, R. Griess Jr. and G. Hoehn,
Framed vertex operator algebras, codes and the moonshine module,
{\em Comm. Math. Phys.} {\bf 193} (1998), 407-448.
\bibitem[DGL]{DGL} C. Dong, R. Griess Jr. and C. Lam, Uniqueness results of the moonshine
vertex operator algebra, {\em American Journal of Math.} {\bf 129}
(2007), 583-609.

\bibitem[DJ]{DJ} C. Dong and C. Jiang, Bimodules associated to  vertex operator algebra,
{\em Math. Z.} {\bf 289} (2008), 799-826.



\bibitem[DL]{DL} C. Dong and J. Lepowsky, Generalized Vertex
Algebras and Relative Vertex Operators, {\em Progress in Math.} Vol. 112,
Birkh\"{a}user, Boston 1993.

\bibitem[DLM1]{DLM1} C. Dong, H. Li and G. Mason,
Regularity of rational vertex operator algebras, {\em  Adv. Math.} {\bf 132} (1997), 148-166.

\bibitem[DLM2]{DLM2} C. Dong, H. Li and G. Mason,
Twisted representations of vertex operator algebras, {\em Math. Ann.}
{\bf  310} (1998), 571--600.

\bibitem[DLM3]{DLM3} C. Dong, H. Li and G. Mason,
Vertex operator algebras and associative algebras, {\em  J. Algebra}
{\bf 206} (1998), 67-96.

\bibitem[DLM4]{DLM4} C. Dong, H. Li and G. Mason, Twisted representations of vertex operator algebras and associative algebras, {\em Int. Math. Res. Not.} {\bf 8} (1998), 389-397.

\bibitem[DMZ]{DMZ} C. Dong, G. Mason and Y. Zhu, Discrete series of the
Virasoro algebra and the moonshine module, {\it Proc. Symp. Pure.
Math., AMS} {\bf 56} II (1994), 295--316.

\bibitem[FFR]{FFR}
Alex J. Feingold, Igor B. Frenkel and John F. X. Ries, Spinor
Construction of Vertex Operator Algebras, Triality, and
$E_{8}^{(1)}$, {\em Contem. Math.} {\bf 121}, 1991.
\bibitem[FF]{FF} B. Feigin and D. Fuchs, Verma Modules over the Virasoro Algebra,
{\em Lect. Notes in Math.} {\bf 1060}, Springer, (1984),  230-245.
\bibitem[FHL]{FHL} I. B. Frenkel, Y. Huang and J. Lepowsky, On
axiomatic approaches to vertex operator algebras and modules,
{\it Memoirs AMS} {\bf 104}, 1993.

\bibitem[FLM]{FLM} I. B. Frenkel, J. Lepowsky and A. Meurman,
Vertex Operator Algebras and the Monster, {\em Pure and Applied
Math.,} Vol. {\bf 134}, Academic Press, 1988.

\bibitem[FZ]{FZ} I. B. Frenkel and Y. Zhu, Vertex operator algebras
associated to  representations of affine and Virasoro algebras, {\em Duke
 Math. J.} {\bf 66} (1992), 123-168.

\bibitem[FQS]{FQS} D. Friedan, Z. Qiu and S. Shenker,
Details of the non-unitarity proof for highest weight representations of
Virasoro Algebra, {\em  Comm. Math. Phys.} {\bf 107} (1986), 535-542.

\bibitem[GKO]{GKO} P. Goddard, A. Kent and D. Olive,
Unitary representations of the Virasoro Algebra and super-Virasoro algebras
{\em Comm. Math. Phys.} {\bf 103} (1986), 105-119.

\bibitem[KR]{KR} V. Kac and A. Raina, Highest Weight Representations of Infinite Dimensional Lie Algebras, Advanced Series in Mathematical Physics, Vol. 2,
World Scientific, 1987.
\bibitem[KW]{KW}
V. Kac and W. Wang, Vertex operator superalgebras and
representations, {\em Contem. Math., AMS} {\bf Vol. 175} (1994), 161-191.

\bibitem[LSY]{LSY}
  C.H. Lam, S. Sakuma and H. Yamauchi, Ising vectors and automorphism
  groups of commutant subalgebras related to root systems,
  {\it Math. Z.} {\bf 255} vol.3 (2007), 597--626.

\bibitem[LS]{LS} C.H. Lam and H. Shimakura, Ising vectors in the vertex operator algebra $V^+_\Lambda$ associated with the Leech lattice $\Lambda$, \emph{Int. Math. Res. Not.}  IMRN 2007, no. 24, Art. ID rnm132, 21 pp.

\bibitem[LYY1]{LYY1}
  C.H. Lam, H. Yamada and H. Yamauchi, McKay's observation and vertex
  operator algebras generated by two conformal vectors of central
  charge 1/2. {\it Int. Math. Res. Papers} {\bf 3} (2005), 117--181.

\bibitem[LYY2]{LYY2} C. H. Lam, H. Yamada and H. Yamauchi, Vertex operator
algebras, extended $E_{8}$ diagram, and McKay's observation on the Monster simple group, \emph{Trans. AMS }\textbf{359} (2007) 4107-4123.

\bibitem[LY]{LY} C. H. Lam and  H. Yamauchi, A characterization of the moonshine vertex operator algebra by means of Virasoro frames, {\em  Int. Math. Res. Not.}  {\bf IMRN 2007} (2007), No.1, Art. ID rnm003, 10 pp

\bibitem[LL]{LL} J. Lepowsky and H. Li, Introduction to Vertex Operator Algebras, and Their Representations, {\em Progress in Math.} Vol. 227,
Birkh\"{a}user, Boston 2004.

\bibitem[L1]{L1} H. Li, Local systems of vertex operators, vertex superalgebras and modules, {\em J. Pure Appl. Algebra} {\bf 109} (1996), 143-195.
\bibitem[L2]{L2}H. Li, Local systems of twisted vertex operators, vertex operator superalgebras and twisted modules, {\em  Contemp. Math. AMS.}
{\bf  193} (1996), 203-236.
\bibitem[L3]{L3} H. Li, The regular representation, Zhu's $A(V)$-theoy and the induced modules, \emph{J. Algebra} \textbf{238} (2001), 159-193.


\bibitem[M1]{M1}  M. Miyamoto, Binary codes and vertex operator (super)algebras,
\emph{J. Algebra} \textbf{181} (1996), 207-222.

\bibitem[M2]{M2} M. Miyamoto, Representation theory of code vertex operator
algebra, {\em  J. Algebra} {\bf  201} (1998), 115--150.

\bibitem[M3]{M3} M. Miyamoto, A new construction of the Moonshine vertex
operator algebra over the real number field, {\em Ann. of Math.}
{\bf 159} (2004), 535--596.


\bibitem[W]{W} W. Wang, Rationality of Virasoro vertex operator algebras,
{\it Int. Math. Res. Not.} {\bf 71} (1993), 197-211.
\bibitem[X]{X} X. Xu, Introduction to vertex operator superalgebras and their modules,
{\em Mathematics and its Applications} {\bf 456}, Kluwer Academic, Dordrecht, 1998.

\bibitem[Z]{Z} Y. Zhu, Modular invariance of characters of vertex operator algebras,
{\em J. Amer, Math. Soc.} {\bf 9} (1996), 237-302.
\end{thebibliography}
\end{document}